\documentclass[preprint,12pt,onecolumn,a4paper]{article}
\usepackage{amsmath}
\usepackage{amssymb}
\usepackage[english]{babel}
\usepackage{tikz}
\usetikzlibrary{intersections, shapes}
\usetikzlibrary{arrows}
\usepackage{subcaption}
\usepackage{booktabs,multirow}
\usepackage{bigdelim}
\usepackage[hidelinks]{hyperref}
\usepackage{array}
\usepackage{pgfplots}
\usepackage[leftcaption]{sidecap}
\sidecaptionvpos{figure}{t}
\usepackage{graphicx,bm,xcolor}

\title{Inexact IETI-DP for conforming isogeometric multi-patch discretizations}
\author{Rainer Schneckenleitner\footnote{\texttt{schneckenleitner@numa.uni-linz.ac.at},
		Institute of Computational Mathematics, Johannes
		Kepler University Linz, Austria}\; and Stefan Takacs\footnote{\texttt{stefan.takacs@ricam.oeaw.ac.at},
		Johann Radon Institute Institute for Computational and Applied Mathematics, Austrian Academy of Sciences, Linz, Austria}}
\date{}

\begin{document}
\maketitle

\begin{abstract}
In this paper, we investigate Dual-Primal Isogeometric Tearing and Interconnecting (IETI-DP) methods for conforming Galerkin discretizations on multi-patch computational domains with inexact subdomain solvers. Recently, the authors have proven a condition number estimate for a IETI-DP method using sparse LU factorizations for the subdomain problems that is explicit, among other parameters, in the grid size and the spline degree. In the present paper, we replace the sparse LU factorizations by fast diagonalization based preconditioners to get a faster IETI-DP method while maintaining the same explicit condition number bound. 
\end{abstract}

\section{Introduction}

The discretization of partial differential equations leads usually to large-scale linear systems.
We are interested in a fast solver for linear systems that are obtained from the discretization of boundary value problems using Isogeometric Analysis (IgA; \cite{HughesCottrellBazilevs:2005}) schemes. We consider computational domains that are composed of multiple non-overlapping patches, for which
FETI-DP type algorithms are a canonical choice. Adaptations of FETI-DP, introduced in~\cite{FarhatLesoinneLeTallecPiersonRixen:2001a}, have already been made to IgA, see, e.g.,~\cite{KleissPechsteinJuttlerTomar:2012,HoferLanger:2016a}. This approach is sometimes called Dual-Primal Isogeometric Tearing and Interconnecting (IETI-DP).
Recently, a convergence analysis for IETI-DP methods that is explicit in the grid sizes, the patch diameters, the spline degree and other parameters like the smoothness of the splines within the patches or the number of patches was carried out for a conforming Galerkin IgA discretization, see~\cite{SchneckenleitnerTakacs:2019}. There, the authors considered a Schur complement IETI-DP method, where the subdomain problems are solved with sparse direct solvers. Large subdomain problems slow down the overall algorithm and require a lot of memory resources. The saddle point formulation of IETI-DP allows the use of inexact local solvers. The successful use of inexact solvers for FETI-DP has already been demonstrated in~\cite{HoferLangerTakacs:2018a,KlawonnRheinbach:2007}. In this paper, we use the fast diagonalization (FD) method introduced in~\cite{SangalliTani:2016} as solver for the subdomain-local subproblems. We show that the inexact IETI-DP version obeys the same condition number bound that also holds for the standard IETI-DP method based on exact solvers for the local problems from~\cite{SchneckenleitnerTakacs:2019}.

The structure of the paper is as follows. Section~\ref{sec:2} is devoted to the introduction of the model problem and the IETI-DP solver. Numerical results are presented in Section~\ref{sec:3}.

\section{Model problem and its solution}
\label{sec:2}
Let $\Omega \subset \mathbb{R}^2$ be an open bounded Lipschitz domain with boundary $\partial \Omega$ and $f\in L_2(\Omega)$ be a given source function. We consider the following model problem: Find $u\in H^1_0(\Omega)$ such that 
\begin{equation}\label{eq:var}
	\int_{\Omega} \nabla u \cdot \nabla v \; \mathrm{d}x = \int_{\Omega} fv \; \mathrm{d}x \qquad \text{for all } v\in H^1_0(\Omega).
\end{equation}
We assume that $\Omega$ is composed of $K$ non-overlapping patches $\Omega^{(k)}$ that are parameterized with geometry mappings 
\[
	G_k : \widehat{\Omega} := (0,1)^d \rightarrow \Omega^{(k)}:=G_k(\widehat{\Omega}),
\]
where the interface between any two patches is either a common vertex or a common edge, cf.~\cite[Ass.~2]{SchneckenleitnerTakacs:2019}. Additionally, we assume that the number of patches sharing a vertex 
is uniformly bounded, cf.~\cite[Ass. 3]{SchneckenleitnerTakacs:2019}. Moreover, we assume that there
is a constant $C_G$ such that $\|\nabla G_k\|_{L_\infty(\widehat\Omega)}\le C_G \, \mbox{diam}(\Omega^{(k)})$
and $\|\nabla G_k^{-1}\|_{L_\infty(\widehat\Omega)}\le C_G \, \mbox{diam}(\Omega^{(k)})^{-1}$,
see~\cite[Ass. 1]{SchneckenleitnerTakacs:2019}, and that we have quasi-uniform grids~\cite[Ass. 4]{SchneckenleitnerTakacs:2019}.
The local discretization spaces on the parameter domain $\widehat{\Omega}$ are tensor-product B-splines spaces, obtained using the Cox-de Boor formula. The local discretization spaces on the physical patches $\Omega^{(k)}$ are obtained by the pull-back principle. We assume that the geometry mappings as well as the discretizations agree on all interfaces between the patches, cf.~\cite[Ass. 5]{SchneckenleitnerTakacs:2019}. So, we are able to set up a fully matching discretization with
the function space
\[
		V = \{ v\in H^1_0(\Omega):v\circ G_k \mbox{ is a B-spline function}\}.
\]
The corresponding discrete problem is obtained by restricting~\eqref{eq:var} to this
space.

In the following, we introduce the IETI-DP solver. The patches from the definition of the
computational domain provide a canonical choice of substructures which we use to set up  the solver.
By assembling the variational problem~\eqref{eq:var} on the individual patches separately, we obtain
independent linear systems
\[
	A^{(k)} \underline u^{(k)} = \underline f^{(k)}
	\quad\mbox{for}\quad
	k=1,\ldots,K,
\]
where $A^{(k)}$ is the local stiffness matrix and $\underline f^{(k)}$ is the
local source vector. We use the vertex values as primal degrees of freedom, see~\cite[Alg.~A]{SchneckenleitnerTakacs:2019}. By eliminating the primal degrees of freedom, we obtain
the matrices $A_{\Delta\Delta}^{(k)}$ and the vectors $\underline f_{\Delta}^{(k)}$ and
$\underline u_{\Delta}^{(k)}$.
Moreover, we introduce jump matrices $B^{(k)}$, where $\sum_{k=1}^K B^{(k)} \underline u^{(k)}=0$ ensures the continuity between the patches (except the continuity at the corners), in the usual way,
see~\cite[Section 3]{SchneckenleitnerTakacs:2019}. The matrices $B_{\Delta}^{(k)}$ are obtained from
$B^{(k)}$ again by eliminating the primal degrees of freedom.

The $A^{(k)}$-orthogonal primal basis representation $\Psi^{(k)}$ is characterized by the linear system
\[
			\begin{pmatrix} A^{(k)} & (C^{(k)})^\top \\   C^{(k)} \end{pmatrix}
		    \begin{pmatrix}  \Psi^{(k)} \\ \Delta^{(k)} \end{pmatrix}
		    =
			\begin{pmatrix} 0 \\ R_c^{(k)} \end{pmatrix},
\]
where $C^{(k)}$ is a matrix that selects the coefficients corresponding to the primal
degrees of freedom and
$R_c^{(k)}$ is a binary local-to-global mapping of the primal constraints.
We compute $\Psi^{(k)}$ using a MINRES solver, preconditioned with
\[
	\widetilde{P}^{(k)} = 
	\begin{pmatrix}
		(\widehat{A}_M^{(k)})^{-1} & \\
		& (C^{(k)} (\widehat{A}_M^{(k)})^{-1} C^{(k)})^{-1}
	\end{pmatrix},
\]
where $\widehat{A}_M^{(k)} := \widehat{A}^{(k)} + \gamma_k \widehat{M}^{(k)}\underline{e}_h^{(k)} (\underline{e}_h^{(k)})^\top \widehat{M}^{(k)}$ is the stiffness matrix on the parameter domain $\widehat{A}^{(k)}$, corrected by a term involving the mass matrix $\widehat{M}^{(k)}$ on the parameter domain and the vector $\underline{e}_h^{(k)}$, which represents the constant function with value $1$. We choose $\gamma_k=1$ if $\Omega^{(k)}$ does not contribute to the Dirichlet boundary $\partial \Omega$ and $\gamma_k=0$ otherwise. The application of $(\widehat{A}_M^{(k)})^{-1}$ is realized with the fast diagonalization (FD) method, see~\cite{SangalliTani:2016}.
After the computation of all $\Psi^{(k)}$, we obtain the global primal basis representation matrix $\Psi$ by canonical mappings. Next, we set up the primal stiffness matrix, the primal source vector and the corresponding jump matrix given by
\[
		A_\Pi := \sum_{k=1}^K (\Psi^{(k)})^\top A^{(k)} \Psi^{(k)},
		\quad \underline{f}_\Pi:=\Psi^\top \underline{f}
		\quad\mbox{and}\quad
		B_\Pi := \sum_{k=1}^K B^{(k)} \Psi^{(k)}.
\]
The overall IETI-DP saddle point system reads as
follows:
\begin{align}\label{eq:matrix problem}
		\begin{pmatrix}
			A_{\Delta\Delta}^{(1)}  &&&	                     &  ({B}_\Delta^{(1)})^\top           \\
			& \ddots &&                                       & \vdots \\
			&& A_{\Delta\Delta}^{(K)}  &	                     &  ({B}_\Delta^{(K)})^\top           \\
			&&&   A_\Pi   & B_\Pi^\top         \\
			B_\Delta^{(1)} &\cdots &B_\Delta^{(K)}          &  B_\Pi              & 0
		\end{pmatrix}	
	\begin{pmatrix}
		\underline{u}_\Delta^{(1)}  \\
		\vdots \\
		\underline{u}_\Delta^{(K)}  \\
		\underline{u}_\Pi  \\
		\underline{\lambda}  \\
	\end{pmatrix}
	=
	\begin{pmatrix}
		\underline{f}_\Delta^{(1)}  \\
		\vdots \\
		\underline{f}_\Delta^{(K)}  \\
		\underline{f}_\Pi  \\
		0  \\
	\end{pmatrix}
	.
\end{align}
We solve~\eqref{eq:matrix problem} with a preconditioned MINRES method using the preconditioner 
\begin{equation}\label{def:saddle point preconditioner}
	\mathcal{P}:=\mbox{diag }
	(
		Q^{(1)} (\widehat{A}_M^{(1)})^{-1} (Q^{(1)})^\top ,\cdots,
		Q^{(K)} (\widehat{A}_M^{(K)})^{-1} (Q^{(K)})^\top,
		A_\Pi^{-1},
		\widehat M_{\mathrm{sD}}
	)	
	,
\end{equation}
where $Q^{(k)}$ is the $A^{(k)}$-orthogonal projection from the local function space into the space $V_\Delta^{(k)}:=\{\underline u^{(k)} : C^{(k)}\underline u^{(k)}=0\}$ of
functions with vanishing primal degrees of freedom and $\widehat{A}_M^{(k)}$ is defined as above and realized using the FD method.
$A_\Pi^{-1}$ is realized by a direct solver.
$\widehat{M}_\mathrm{sD}$ is an inexact scaled Dirichlet preconditioner, defined by 
\begin{equation}\label{eq:cG:approximate scaled Dirichlet 1}
	\widehat{M}_{\mathrm{sD}} := B_\Gamma {D}^{-1} \widehat{S} {D}^{-1} B_\Gamma^\top,
\end{equation}
where
$\widehat{S} := \mbox{diag}(\widehat{S}^{(1)}, \dots, \widehat{S}^{(K)})$
with
$\widehat S^{(k)} := \widehat A_{\Gamma \Gamma}^{(k)} - \widehat A_{\Gamma \mathrm{I}}^{(k)} (\widehat A_{\mathrm{II}}^{(k)})^{-1} \widehat A_{\mathrm{I} \Gamma}^{(k)}$.
The diagonal matrix $D$ is based on the principle of multiplicity scaling, cf.~\cite{SchneckenleitnerTakacs:2019}. 
The index $\Gamma$ refers to the rows/columns of
$\widehat{A}^{(k)}$ and the columns of $B^{(k)}$ that refer to basis functions
with non-vanishing trace, the index $\mathrm{I}$ refers to the remaining
rows/columns.

Under the presented assumptions, the condition number of the system~\eqref{eq:matrix problem} preconditioned with~\eqref{def:saddle point preconditioner} is bounded by
\[
	C\, p \; 
	\left(1+\log p+\max_{k=1,\ldots,K} \log\frac{H_{k}}{h_{k}}\right)^2,
\]
where $p$ is the spline degree, $H_k$ is the diameter of $\Omega^{(k)}$,
and $h_k$ the grid size on $\Omega^{(k)}$, and the constant $C$ only depends on the constant
$C_G$, the quasi-uniformity constant and the maximum number of patches sharing a vertex, see~\cite{SchneckenleitnerTakacs:2019}.

\section{Numerical results}\label{sec:3}
\begin{figure}[h]
	\centering
	\includegraphics[height=3.5cm]{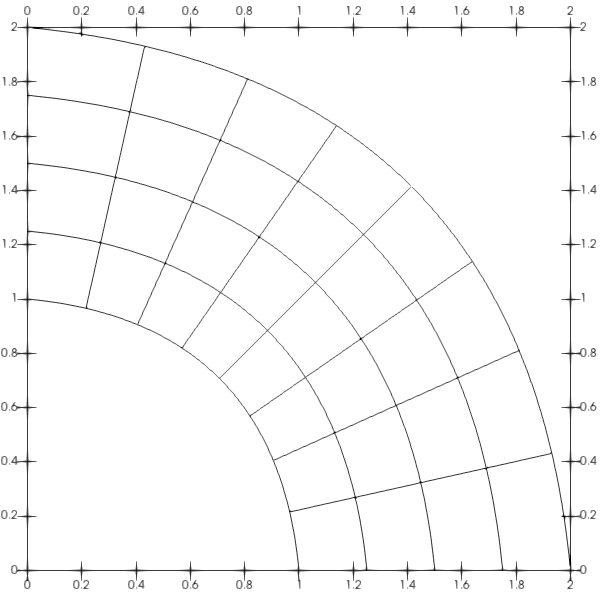}
	\qquad\qquad
	\includegraphics[height=3.5cm]{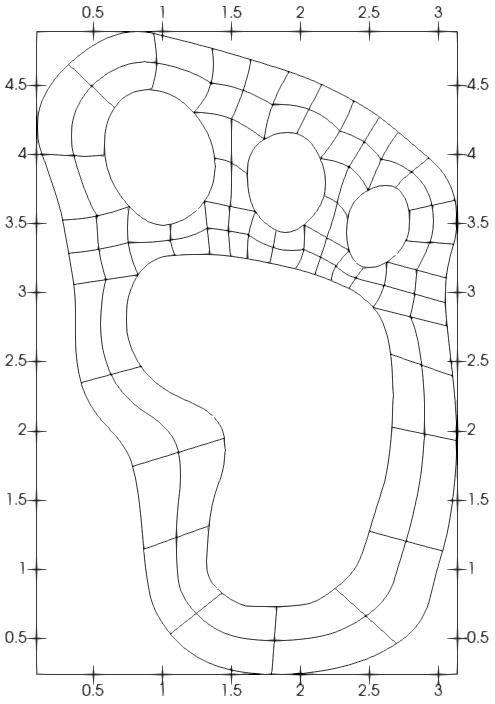}
	\caption{Computational domains; Quarter annulus (left); Yeti-footprint (right)}
	\label{fig:domain}
\end{figure}

In this section, we show numerical results of the proposed inexact IETI-DP method on the computational domains as given in Fig.~\ref{fig:domain}. The first domain is a quarter annulus consisting of $32$ patches and the second one is the Yeti-footprint decomposed into $84$ patches. 

We consider the model problem
\[
- \Delta u(x,y)  = 2\pi^2 \sin(\pi x)\sin(\pi y)
\quad \mbox{for}\quad
(x,y)\in\Omega
\]
with homogeneous Dirichlet boundary conditions.
Within each patch, we use B-splines of degree $p$ and maximum smoothness $C^{p-1}$. The coarsest discretization space ($r=0$) is the space of patchwise global polynomials, only the more
rectangular patches of the Yeti-footprint have one inner knot on each of the longer sides of
the patches. The subsequent refinements $r=1,2,\dots$ are obtained via uniform refinement steps.

All numerical experiments have been carried out using
the C++ library G+Smo\footnote{\url{https://github.com/gismo/gismo}}, the CPU times have been
recorded on the Radon1\footnote{\url{https://www.ricam.oeaw.ac.at/hpc/}} Cluster in Linz. 

\begin{figure}[ht]
	\begin{center}
		\resizebox{5cm}{5cm}{%
			\begin{tikzpicture}
				\begin{axis}[
					xlabel={Refinement level $r$},
					ylabel={Solving times [sec]},
					xmin=5, xmax=9,
					ymin=1, ymax=1500,
					xtick={5,6,7,8,9},
					ytick={1,2,5,10,20,50,100,200,500,1000},
					legend pos=north west,
					ymajorgrids=true,
					grid style=dashed,
					legend columns=3,
					legend pos=north west,
					legend image post style={only marks},
					log ticks with fixed point,
					ymode=log,
					width=6cm
					]
					
					\addplot[color=blue,
					mark=triangle]
					coordinates {
						(5,1.2)(6,4.2)(7,25.0)(8,126.3)(9,708)
					};
					\addlegendentry{{\tiny MFD}};
					
					\addplot[color=red,
					mark=square]
					coordinates {
						(5,6.7)(6,43.0)(7,216.0)(8,1015.0)
					};
					\addlegendentry{{\tiny MLU}};
					
					\addplot[color=green,
					mark=x]
					coordinates {
						(5,2.7)(6,17.0)(7,81.0)(8,412.0)
					};
					\addlegendentry{{\tiny CGLU}};
					
				\end{axis}
			\end{tikzpicture}
		}
		\quad
		\resizebox{5cm}{5cm}{%
			\begin{tikzpicture}
				\begin{axis}[
					xlabel={Polynomial degree level $p$},
					ylabel={Solving times [sec]},
					xmin=2, xmax=8,
					ymin=10, ymax=500,
					xtick={2,3,4,5,6,7,8},
					ytick={10,20,50,100,200,500},
					legend pos=north west,
					ymajorgrids=true,
					grid style=dashed,
					legend columns=3,
					legend pos=north west,
					legend image post style={only marks},
					log ticks with fixed point,
					ymode=log,
					width=6cm
					]
					
					\addplot[color=blue,
					mark=triangle]
					coordinates {
						(2,14.0)(3,17.0)(4,21.0)(5,25.0)(6,30.0)(7,35.0)(8,41.0)
					};
					\addlegendentry{{\tiny MFD}};
					
					\addplot[color=red,
					mark=square]
					coordinates {
						(2,83.0)(3,129.0)(4,174.0)(5,216.0)(6,248.0)(7,297.0)(8,343.0)
					};
					\addlegendentry{{\tiny MLU}};
					
					\addplot[color=green ,
					mark=x]
					coordinates {
						(2,32.0)(3,51.0)(4,65.7)(5,81.0)(6,99.0)(7,104.0)(8,134.0)
					};
					\addlegendentry{{\tiny CGLU}};
					
				\end{axis}
			\end{tikzpicture}
		}
	\end{center}
	\caption{Solving times for $p=5$ (left) and $r=7$ (right); MFD (blue lines); MLU (red lines); CGLU (green lines); Annulus
		\label{fig:1}}
\end{figure}
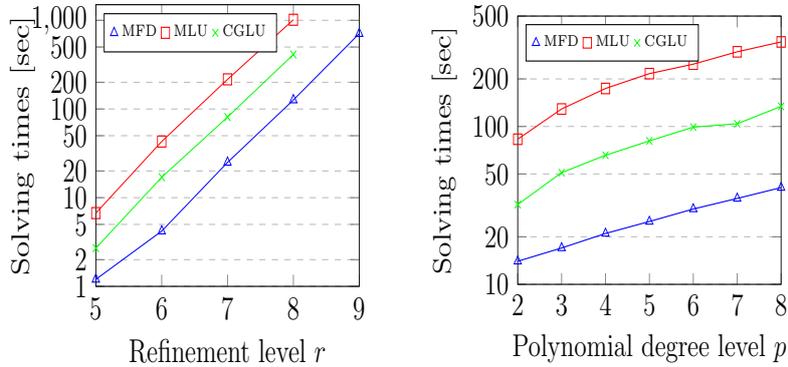

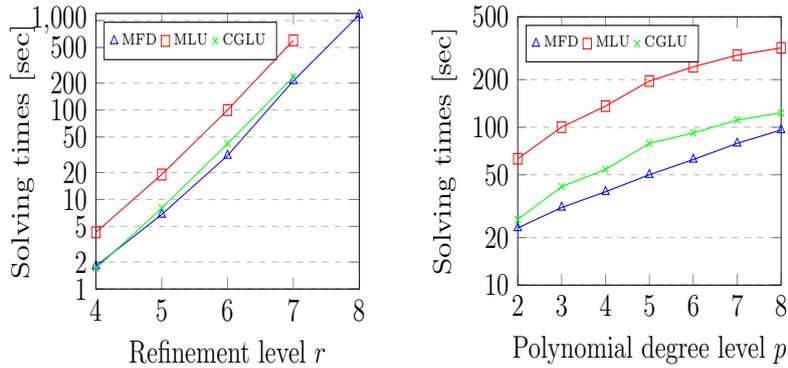
\begin{figure}[ht]
\begin{center}
\resizebox{5cm}{5cm}{%
\begin{tikzpicture}
\begin{axis}[
	xlabel={Refinement level $r$},
	ylabel={Solving times [sec]},
	xmin=4, xmax=8,
	ymin=1, ymax=1200,
	xtick={4,5,6,7,8},
	ytick={1,2,5,10,20,50,100,200,500,1000,1500},
	legend pos=north west,
	ymajorgrids=true,
	grid style=dashed,
	legend columns=3,
	legend pos=north west,
	legend image post style={only marks},
	log ticks with fixed point,
	ymode=log,
	width=6cm
	]
	
	\addplot[color=blue,
	mark=triangle]
	coordinates {
	(4,1.8)(5,6.8)(6,31.0)(7,212.0)(8,1158)
	};
	\addlegendentry{{\tiny MFD}};
	
	\addplot[color=red,
	mark=square]
	coordinates {
	(4,4.3)(5,19.0)(6,100.0)(7,600.5)
	};
	\addlegendentry{{\tiny MLU}};
	
	\addplot[color=green,
	mark=x]
	coordinates {
	(4,1.7)(5,8.0)(6,42.0)(7,234.0)
	};
	\addlegendentry{{\tiny CGLU}};

\end{axis}
\end{tikzpicture}
}
\quad
\resizebox{5cm}{5cm}{%
\begin{tikzpicture}
\begin{axis}[
	xlabel={Polynomial degree level $p$},
	ylabel={Solving times [sec]},
	xmin=2, xmax=8,
	ymin=10, ymax=500,
	xtick={2,3,4,5,6,7,8},
	ytick={1,2,5,10,20,50,100,200,500},
	legend pos=north west,
	ymajorgrids=true,
	grid style=dashed,
	legend columns=3,
	legend pos=north west,
	legend image post style={only marks},
	log ticks with fixed point,
	ymode=log,
	width=6cm
	]
	
	\addplot[color=blue,
	mark=triangle]
	coordinates {
	(2,23.0)(3,31.0)(4,39.0)(5,50.0)(6,62.4)(7,79.0)(8,96.0)
	};
	\addlegendentry{{\tiny MFD}};

	\addplot[color=red,
	mark=square]
	coordinates {
	(2,63.1)(3,100.0)(4,136.0)(5,196.0)(6,241.2)(7,286.0)(8,318.0)
	};
	\addlegendentry{{\tiny MLU}};

	\addplot[color=green ,
	mark=x]
	coordinates {
	(2,26.0)(3,42.0)(4,54.0)(5,79.0)(6,92.0)(7,111.0)(8,123.2)
	};
	\addlegendentry{{\tiny CGLU}};

\end{axis}
\end{tikzpicture}
}
\end{center}
	\caption{Solving times for $p=3$ (left) and $r=6$ (right); MFD (blue lines); MLU (red lines); CGLU (green lines); Yeti-footprint
		\label{fig:2}}
\end{figure}


We compare three different IETI-DP solvers: the proposed solver as introduced in Section~\ref{sec:2} (=MFD), a IETI-DP solver for the saddle point system~\eqref{eq:matrix problem} without the primal degrees of freedom eliminated with
direct solvers for the local subproblems (=MLU), and the Schur complement based approach as introduced in~\cite{SchneckenleitnerTakacs:2019} (=CGLU). We use MINRES as outer solver in the cases MFD and MLU
and conjugate gradient as outer solver for the case CGLU. For MLU and CGLU, we use sparse direct LU solvers from the Pardiso project\footnote{\url{https://www.pardiso-project.org/}} for the local subproblems and for computing the bases $\Psi^{(k)}$. For MFD, the primal bases $\Psi^{(k)}$ are solved with MINRES up to an accuracy of $10^{-8}$. We present the times required for computing a primal basis $\Psi$ which we indicate in the tables with the same letter, the accumulated setup and application times of the different local preconditioners for all patches $K$, abbreviated by $\Theta_S$ and $\Theta_A$, respectively as well as the solving times and
the number of iterations (it.) required by the main Krylov space solver. We start all numerical experiments with zero initial guess and stop the iterations if the ${\ell_2}$-norm of the residual vector is reduced by a factor of $10^{-6}$ compared to the ${\ell_2}$-norm of the right-hand side vector. 

\begin{table}[th]
	\scriptsize
	\newcolumntype{L}[1]{>{\raggedleft\arraybackslash\hspace{-1em}}m{#1}}
	\centering
	\renewcommand{\arraystretch}{1.25}
	\begin{tabular}{r|c|L{3.em}|L{3.0em}|L{3.em}|L{3.em}|L{3.em}|L{2.em}}
		\toprule
		\multicolumn{1}{c|}{\;\footnotesize $ $ }
		&\multicolumn{1}{c|}{\;\footnotesize $r$ \;}
		& \multicolumn{1}{c|}{\footnotesize $\Psi$}
		& \multicolumn{1}{c|}{\footnotesize $\Theta_{S}$} 
		& \multicolumn{1}{c|}{\footnotesize $\Theta_{A}$}
		& \multicolumn{1}{c|}{\footnotesize solving} 
		& \multicolumn{1}{c|}{\footnotesize total}
		& \multicolumn{1}{c}{\footnotesize it.}\\
		\midrule
		\midrule
		MFD & $6$   
		& $1.2$ & $0.1$ & $0.4$ & $4.2$ & $5.5$ & $71$  \\
		MLU &  $ $ 
		& $0.8$ & $7.6$ & $21.2$ & $43.0$ & $51.4$ & $37$  \\
		CGLU &  $ $ 
		& $0.8$ & $7.6$ & $9.5$ & $17.0$ & $15.4$ & $15$ \\
		\midrule
		MFD & $7$   
		& $8.0$ & $0.3$ & $4.9$ & $25.0$ & $33.3$ & $80$ \\
		MLU &  $ $
		& $4.0$ & $42.3$ & $ 106.4$ & $216.0$ & $262.3$ & $39$ \\
		CGLU &  $ $ 
		& $4.0$ & $42.3$ & $45.4$ & $81.0$ & $127.3$ & $15$  \\
		\midrule
		MFD & $8$   
		& $35.5$ & $2.0$ & $26.9$ & $126.3$ & $163.8$ & $88$\\
		MLU &  $ $ 
		& $18.1$ & $243.6$ & $503.8$ & $1015.0$ & $1276.7$ & $41$ \\
		CGLU &  $ $ 
		& $18.0$ & $242.1$ & $225.3$ & $412.0$ & $672.1$ & $17$\\
		\bottomrule
	\end{tabular}
	\captionof{table}{Alg.~A; $p = 5$; timings in sec.; Annulus
		\label{tab:Annulus:Alg. A:p = 5}}
\end{table}

\begin{table}[th]
	\scriptsize
	\newcolumntype{L}[1]{>{\raggedleft\arraybackslash\hspace{-1em}}m{#1}}
	\centering
	\renewcommand{\arraystretch}{1.25}
	\begin{tabular}{r|c|L{3.em}|L{3.0em}|L{3.em}|L{3.em}|L{3.em}|L{2.em}}
		\toprule
		\multicolumn{1}{c|}{\;\footnotesize $ $ }
		&\multicolumn{1}{c|}{\;\footnotesize $r$ \;}
		& \multicolumn{1}{c|}{\footnotesize $\Psi$}
		& \multicolumn{1}{c|}{\footnotesize $\Theta_{S}$} 
		& \multicolumn{1}{c|}{\footnotesize $\Theta_{A}$}
		& \multicolumn{1}{c|}{\footnotesize solving} 
		& \multicolumn{1}{c|}{\footnotesize total}
		& \multicolumn{1}{c}{\footnotesize it.}\\
		\midrule
		\midrule
		MFD & $6$   
		& $2.5$ & $0.2$ & $0.5$ & $8.2$ & $10.9$ & $76$  \\
		MLU &  $ $ 
		& $1.2$ & $17.9$ & $33.8$ & $72.0$ & $91.1$ & $39$  \\
		CGLU &  $ $ 
		& $1.2$ & $17.8$ & $14.4$ & $26.0$ & $45.0$ & $15$ \\
		\midrule
		MFD & $7$   
		& $15.1$ & $0.4$ & $5.4$ & $41.0$ & $56.5$ & $84$ \\
		MLU &  $ $
		& $5.8$ & $105.0$ & $163.3$ & $343.0$ & $453.8$ & $41$ \\
		CGLU &  $ $ 
		& $5.8$ & $105.0$ & $73.8$ & $134.0$ & $244.8$ & $17$  \\
		\midrule
		MFD & $8$   
		& $60.1$ & $2.3$ & $29.9$ & $200.0$ & $262.4$ & $94$\\
		MLU &  $ $ 
		& \multicolumn{6}{c}{OoM} \\
		CGLU &  $ $ 
		& \multicolumn{6}{c}{OoM}\\
		\bottomrule
	\end{tabular}
	\captionof{table}{Alg.~A; $p = 8$; timings in sec.; Annulus
		\label{tab:Annulus:Alg. A:p = 8}}
\end{table}

The Tables~\ref{tab:Annulus:Alg. A:p = 5} and~\ref{tab:Annulus:Alg. A:p = 8} present the timings of the algorithms on the quarter annulus domain. We observe that the solving and total times required by MFD are three to five times smaller compared to the other methods MLU and CGLU. MFD is much faster than the other methods despite the fact that the required number of iterations are up to approximately six times larger. One disadvantage of MFD is the computation of the primal basis $\Psi$. The tables show a larger computation time when using MFD. This is a weakness of the classical preconditioned MINRES method when applied to problems with multiple right-hand sides. In general, we have to solve systems with four right-hand sides to compute $\Psi^{(k)}$. Moreover, we see that the setup and application of the FD based preconditioner is much faster compared to the factorizations of 
\[
	\begin{pmatrix}
		A^{(k)} & (C^{(k)})^\top \\
		C^{(k)} & 
	\end{pmatrix}
\]
and their application. Table~\ref{tab:Annulus:Alg. A:p = 8} shows another advantage of the MFD method. Since
its memory footprint is smaller, it also provides a solution vector to the considered linear system for the refinement level $r = 8$.  
The plots in Fig.~\ref{fig:1} visualize solving times of the IETI-DP solvers on the quarter annulus domain. We mark the performance of MFD with blue lines and triangles, MLU with red lines and squares and the performance of CGLU is indicated with green lines and crosses. In both graphs, we observe that MFD is the fastest algorithm. In the left plot, we see that for spline degree $p=5$, the solving times increase rather linearly with respect to the number of unknowns. Moreover, the left graph shows that MFD computes the solution for the linear system even for refinement level $r=9$ $(\approx 8.5\mbox{M} \, \mbox{dofs})$. 
In the right graph, we present the solving times with respect to the spline degree for refinement level $r=7$. The solving times for the three IETI-DP solvers increase about linearly with the spline degree.

\begin{table}[tb]
	\scriptsize
	\newcolumntype{L}[1]{>{\raggedleft\arraybackslash\hspace{-1em}}m{#1}}
	\centering
	\renewcommand{\arraystretch}{1.25}
	\begin{tabular}{r|c|L{3.em}|L{3.0em}|L{3.em}|L{3.em}|L{3.em}|L{2.em}}
		\toprule
		\multicolumn{1}{c|}{\;\footnotesize $ $ }
		&\multicolumn{1}{c|}{\;\footnotesize $r$ \;}
		& \multicolumn{1}{c|}{\footnotesize $\Psi$}
		& \multicolumn{1}{c|}{\footnotesize $\Theta_{S}$} 
		& \multicolumn{1}{c|}{\footnotesize $\Theta_{A}$}
		& \multicolumn{1}{c|}{\footnotesize solving} 
		& \multicolumn{1}{c|}{\footnotesize total}
		& \multicolumn{1}{c}{\footnotesize it.}\\
		\midrule
		\midrule
		MFD & $5$   
		& $0.7$ & $0.1$ & $0.8$ & $6.8$ & $7.6$ & $213$  \\
		MLU &  $ $ 
		& $0.3$ & $2.3$ & $9.4$ & $19.0$ & $21.6$ & $45$  \\
		CGLU &  $ $ 
		& $0.3$ & $2.3$ & $4.5$ & $8.0$ & $10.6$ & $20$ \\
		\midrule
		MFD & $6$   
		& $3.4$ & $0.2$ & $5.0$ & $31.0$ & $34.6$ & $242$ \\
		MLU &  $ $
		& $1.2$ & $10.5$ & $49.6$ & $100.0$ & $111.7$ & $51$ \\
		CGLU &  $ $ 
		& $1.2$ & $10.4$ & $23.1$ & $42.0$ & $53.6$ & $22$  \\
		\midrule
		MFD & $7$   
		& $21.8$ & $1.2$ & $ 52.7$ & $212.0$ & $235.0$ & $274$\\
		MLU &  $ $ 
		& $6.8$ & $50.3$ & $295.2$ & $600.5$ & $657.6$ & $55$\\
		CGLU &  $ $ 
		& $6.8$ & $50.1$ & $126.1$ & $234.0$ & $290.9$ & $22$\\
		\bottomrule
	\end{tabular}
	{\captionof{table}{Alg.~A; $p = 3$; timings in sec.; Yeti-footprint}
	\label{tab:Vertex, p = 3}}
\end{table}

\begin{table}[tb]
	\scriptsize
	\newcolumntype{L}[1]{>{\raggedleft\arraybackslash\hspace{-1em}}m{#1}}
	\centering
	\renewcommand{\arraystretch}{1.25}
	\begin{tabular}{r|c|L{3.em}|L{3.0em}|L{3.em}|L{3.em}|L{3.em}|L{2.em}}
		\toprule
		\multicolumn{1}{c|}{\;\footnotesize $ $ }
		&\multicolumn{1}{c|}{\;\footnotesize $r$ \;}
		& \multicolumn{1}{c|}{\footnotesize $\Psi$}
		& \multicolumn{1}{c|}{\footnotesize $\Theta_{S}$} 
		& \multicolumn{1}{c|}{\footnotesize $\Theta_{A}$}
		& \multicolumn{1}{c|}{\footnotesize solving} 
		& \multicolumn{1}{c|}{\footnotesize total}
		& \multicolumn{1}{c}{\footnotesize it.}\\
		\midrule
		\midrule
		MFD & $5$   
		& $2.2$ & $0.2$ & $1.2$ & $20.0$ & $22.4$ & $249$  \\
		MLU &  $ $ 
		& $0.6$ & $10.5$ & $23.5$ & $48.0$ & $59.1$ & $51$  \\
		CGLU &  $ $ 
		& $0.6$ & $10.4$ & $10.9$ & $19.0$ & $30.0$ & $22$ \\
		\midrule
		MFD & $6$   
		& $9.5$ & $0.4$ & $5.9$ & $79.0$ & $88.9$ & $282$ \\
		MLU &  $ $
		& $3.0$ & $44.8$ & $140.5$ & $286.0$ & $333.8$ & $57$ \\
		CGLU &  $ $ 
		& $3.0$ & $44.8$ & $60.6$ & $111.0$ & $158.8$ & $23$  \\
		\midrule
		MFD & $7$   
		& $54.2$ & $1.5$ & $63.2$ & $414.0$ & $469.7$ & $309$\\
		MLU &  $ $ 
		& $14.1$ & $265.3$ & $739.1$ & $1521.0$ & $1800.4$ & $61$\\
		CGLU &  $ $ 
		& $13.8$ & $261.7$ & $304.4$ & $570.0$ & $845.5$ & $25$\\
		\bottomrule
	\end{tabular}
	{\captionof{table}{Alg.~A; $p = 7$; timings in sec.; Yeti-footprint}
	\label{tab:Vertex, p = 7}}
\end{table}

The plots in Fig.~\ref{fig:2} show solving times of the IETI-DP solvers on the Yeti-footprint. We marked the performance of the different IETI-DP solvers as above in the experiments on the quarter annulus. The plot on the left shows the increase of the solving time with respect to the refinement level with polynomial degree $p = 3$ and the plot on the right shows the increase of the solving time with respect to the polynomial degree, where we have fixed the refinement level to $r=6$. As on the quarter annulus, we see that MFD is superior compared to MLU and CGLU also on the Yeti-footprint with respect to the solving times and the smaller memory footprint of MFD allows us to consider larger problems. In both plots, we observe similar growth rates of the solving time for all three solvers as in Fig.~\ref{fig:1}.
In the Tables~\ref{tab:Vertex, p = 3} and~\ref{tab:Vertex, p = 7}, we present and compare the required timings for the polynomial degrees $p = 3$ and $p = 7$ for different refinement levels on the Yeti-footprint.

To conclude, we presented a fast IETI-DP method which allows the incorporation of inexact solvers for the local subproblems while maintaining the condition number bound as established in~\cite{SchneckenleitnerTakacs:2019}. It is beneficial both because of its smaller memory footprint
and its faster convergence for the model problems.
 
\textbf{Acknowledgments.} The first author was supported by the Austrian Science Fund (FWF): S117-03 and W1214-04. Also, the second author has received support from the Austrian Science	Fund (FWF): P31048.

\bibliography{literature}

\end{document}